\documentclass[11pt]{amsart}

\usepackage{amssymb,amsfonts,amscd,amsbsy}
\usepackage[mathscr]{eucal}
\usepackage{url}

\theoremstyle{plain}
\newtheorem{thm}{Theorem}[section]
\newtheorem{lem}[thm]{Lemma}
\newtheorem{prop}[thm]{Proposition}
\newtheorem{cor}[thm]{Corollary}
\theoremstyle{definition}

\newtheorem{rem}[thm]{Remark}

\newcommand{\SigmaP}{{\sf Sigma}}

\newcounter{mmacnt}
\def\restartmma{\setcounter{mmacnt}{0}}
\restartmma \catcode`|=\active
\def|#1|{\mathrm{#1}}
\catcode`|=12
\newenvironment{mma}{
 \par\smallskip
 \catcode`|=\active
 \parskip=0pt\parindent=0pt % locally
 \small
 \def\In##1\\{%
   \def\linebreak{\hfill\break\null\qquad}%
   \refstepcounter{mmacnt}
   \hangindent=2.5em\hangafter=0
   \leavevmode
   \llap{\tiny\sffamily In[\arabic{mmacnt}]:=\kern.5em}%
   \mathversion{bold}\footnotesize$\tt\bf\displaystyle##1$\normalsize
   \mathversion{normal}\par
 }%
 \def\Print##1\\{%
   \def\linebreak{\hfill\break}%
   \hangindent=2.5em\hangafter=0
   \leavevmode \footnotesize##1\par}%
 \def\Out##1\\{%
   \def\linebreak{$\hfill\break\null\hfill$}%
   \kern\abovedisplayskip\par
   \hangindent=2.5em\hangafter=0
   \leavevmode
   \llap{\tiny\sffamily Out[\arabic{mmacnt}]=\kern.5em}
   \footnotesize$\displaystyle\tt##1$\normalsize\hfill\null\par
   \kern\belowdisplayskip
 }%
 \def\Warning##1##2\\{%
   \def\linebreak{\hfill\break}%
   \hangindent=2.5em\hangafter=0
   \leavevmode
   {\scriptsize##1 : ##2}\par}%
}{%
 \par\smallskip
}

\newcommand{\myOut}[1]{{\sffamily Out[#1]}}

\def\MLabel#1{{\refstepcounter{mmacnt}\label{#1}}\addtocounter{mmacnt}{-1}}

\begin{document}
\title[Gaussian Hypergeometric series]{Gaussian Hypergeometric series and supercongruences}
\author{Robert Osburn and Carsten Schneider}

\address{School of Mathematical Sciences, University College Dublin, Belfield, Dublin 4, Ireland}

\address{IH{\'E}S, Le Bois-Marie, 35, route de Chartres, F-91440 Bures-sur-Yvette, FRANCE}

\address{Research Institute for Symbolic Computation, J. Kepler University Linz, Altenberger Str. 69, A-4040 Linz, Austria}

\email{robert.osburn@ucd.ie, osburn@ihes.fr}

\email{Carsten.Schneider@risc.uni-linz.ac.at}

\thanks{The second author was supported by the SFB-grant F1305 and the grant P16613-N12
of the Austrian FWF}


\date{November 8, 2007}

\begin{abstract}  Let $p$ be an odd prime. In 1984, Greene introduced the notion of hypergeometric functions over finite fields. Special values of these functions have been of interest as they are related to the number of $\mathbb{F}_{p}$ points on algebraic varieties and to Fourier coefficients of modular forms. In this paper, we explicitly determine these functions modulo higher powers of $p$ and discuss an application to supercongruences. This application  uses two non-trivial generalized Harmonic sum identities discovered using the computer summation package ~\SigmaP. We illustrate the usage of ~\SigmaP~ in the discovery and proof of these two identities.

\end{abstract}

\maketitle

\section{Introduction}

In \cite{g1} and \cite{g2}, Greene defined general hypergeometric series over finite fields. His aim was to show that these functions satisfy properties analogous to classical hypergeometric series. For example, the four major evaluations of the ordinary hypergeometric series ${}_{3}F_{2}$ due to Saalsch{\"u}tz, Dixon, Watson, and Whipple \cite{aar} all have finite field interpretations (see page 126 of \cite{g1}). Greene's work was in part motivated by the analogy between Gauss sums and Gamma functions \cite{e}, \cite{k1}, \cite{y}.

His approach has proven to be a powerful technique for character sum evaluations.
Recently, several authors have shown that special values of these functions are related to
the number of points over $\mathbb{F}_{p}$, $p$ an odd prime, of Calabi-Yau threefolds \cite{AO}, traces of Hecke operators \cite{fop}, formulas for Ramanujan's $\tau$-function \cite{p}, and the number of points on a family of elliptic curves \cite{jf}. We should also mention that hypergeometric series over arbitrary fields has been developed \cite{gg}, \cite{ggr}, but their application to number theory has yet to be
investigated.

The purpose of this paper is to further study arithmetic properties of hypergeometric functions over finite fields. In particular, we explicitly determine these functions modulo higher powers of $p$ and then briefly discuss extensions of supercongruences.

We first recall some definitions. Let $\mathbb{F}_{p}$ denote the finite field with $p$ elements. We extend all characters $\chi$ of $\mathbb{F}_{p}^{*}$ to $\mathbb{F}_{p}$ by setting $\chi(0):=0$. Following \cite{g1} and \cite{g2}, we give two definitions. The first definition is the finite field analogue of the binomial coefficient. If $A$ and $B$ are characters of $\mathbb{F}_{p}$, then

\begin{equation} \label{def1}
\binom{A}{B} := \frac{B(-1)}{p}J(A, \bar{B})=\frac{B(-1)}{p}
\sum_{x\in \mathbb{F}_{p}}A(x)\bar{B}(1-x),
\end{equation}

\noindent where $J(\chi, \psi)$ denotes the Jacobi sum if $\chi$ and $\psi$ are characters of
$\mathbb{F}_{p}$. The second definition is the finite field analogue of ordinary hypergeometric functions. If $A_{0}$, $A_{1}$, $\dotso$, $A_{n}$, and $B_{1}$, $\dotso$, $B_{n}$ are
characters of $\mathbb{F}_{p}$, then the Gaussian hypergeometric function over $\mathbb{F}_{p}$ is defined by

\begin{equation} \label{def2}
{}_{n+1}F_{n}
\left(\begin{matrix} A_0, &  A_1, & \dotsc,& A_n \\ \ &   B_1, & \dots,& B_n
\end{matrix}
\mid x\right)_p:=\frac p{p-1}\sum_\chi{A_0\chi \choose\chi}{A_1\chi\choose
B_1\chi}\dots {A_n\chi\choose B_n\chi}\chi(x),
\end{equation}

\noindent where the summation is over all characters $\chi$ of $\mathbb{F}_{p}$. In this paper, we restrict our attention to the case $A_{i}=\phi_{p}$ for all $i$ and $B_{j}=\epsilon_{p}$ for all $j$ where $\phi_{p}$ is the quadratic character and $\epsilon_{p}$ is the trivial character mod $p$. We shall denote this value by ${}_{n+1}F_{n}(\lambda)$. By \cite{g1} and \cite{g2}, $p^{n} {}_{n+1}F_{n}(\lambda) \in \mathbb{Z}$. Before stating the main result, we recall that
for $i$, $n \in \mathbb{N}$, generalized Harmonic sums $H_{n}^{(i)}$ are defined by

$$
H_{n}^{(i)} := \sum_{j=1}^{n} \frac{1}{j^i}
$$

\noindent and $H_{0}^{(i)}:=0$. We now define the quantities

\begin{equation} \label{x}
\begin{aligned}
X(p, \lambda, n) & := \phi_{p}(\lambda) \sum_{j=0}^{\tfrac{p-1}{2}} \binom{\frac{p-1}{2} + j}{j}^{l} \binom{\frac{p-1}{2}}{j}^{l} (-1)^{jl}   \lambda^{-j} \Biggl(1+ 2(n+1)j \Bigl(H_{\frac{p-1}{2} + j}^{(1)} - H_{j}^{(1)} \Bigr) \\
&+ j^{2} \Bigl( \tfrac{n+1}{2}(1+n) \Bigl ( H_{\frac{p-1}{2} + j}^{(1)} - H_{j}^{(1)} \Bigr)^2 - (\tfrac{n+1}{2})\Bigl( H_{\frac{p-1}{2} + j}^{(2)} - H_{j}^{(2)} \Bigr) \Bigr)  \Biggr),
\end{aligned}
\end{equation}

\noindent

\begin{equation} \label{y}
\begin{aligned}
Y(p, \lambda, n) &:= \phi_{p}(\lambda)  \sum_{j=0}^{\tfrac{p-1}{2}} \binom{\frac{p-1}{2} + j}{j}^{l} \binom{\frac{p-1}{2}}{j}^{l} (-1)^{jl}  \lambda^{-jp} \Biggl(1 + (n+1)j \Bigl(H_{\frac{p-1}{2} + j}^{(1)} \\
& - H_{j}^{(1)} \Bigr) -  (\tfrac{n+1}{2})j \Bigl(H_{\frac{p-1}{2} + j}^{(1)} - H_{\frac{p-1}{2} - j}^{(1)} \Bigr) \Biggr),
\end{aligned}
\end{equation}

\noindent

\begin{equation} \label{z}
Z(p,\lambda, n):= \phi_{p}(\lambda) \sum_{j=0}^{\tfrac{p-1}{2}} \binom{2j}{j}^{2l} 16^{-jl}  \lambda^{-jp^2},
\end{equation}

\noindent and

\begin{equation} \label{d}
D(p, \lambda):= \sum_{j=0}^{\frac{p-5}{2}} \frac{j!^2}{\prod_{i=0}^{j} (i+ \frac{1}{2})^2} (j+1)^{2} \lambda^{-j-1},
\end{equation}

\noindent where $l:=\frac{n+1}{2}$. The main result of this paper is the following.

\begin{thm} \label{mains} If $n \geq 2$, then

\begin{equation} \label{n3}
-p^{n} {}_{n+1}F_{n}(\lambda) \equiv (-\phi_{p}(-1))^{n+1} \Bigl[ p^2 X(p,\lambda, n) + pY(p,\lambda, n) + Z(p,\lambda, n) \Bigr] \pmod {p^3}
\end{equation}

\noindent and if $n=1$, then

\begin{equation} \label{n1}
-p {}_{2}F_{1}(\lambda) \equiv p^2 \Bigl[ X(p,\lambda, 1) + D(p, \lambda) \Bigr] + pY(p,\lambda, 1) + Z(p,\lambda, 1) \pmod {p^3}.
\end{equation}

\end{thm}

We note that Theorem \ref{mains} generalizes both Theorem 1 in \cite{a}, where the case $n=2$ was handled modulo $p^2$, and Theorem 2.4 in \cite{lr}. As an application of Theorem \ref{mains}, we prove a supercongruence for the Legendre symbol $\bigl(\frac{-1}{p}\bigr) $. This result generalizes Theorem 1 in \cite{mort}.

\begin{cor} \label{done} Let $p$ be an odd prime.  Then

\begin{equation} \label{leg}
\sum_{n=0}^{\frac{p-1}{2}} \binom{2n}{n}^2 16^{-n} +  \frac{3}{8}p (-1)^{\tfrac{p-1}{2}} \sum_{i=1}^{\frac{p-1}{2}}  \binom{2i}{i} \frac{1}{i}
\equiv  \Biggl(\frac{-1}{p}\Biggr) \pmod {p^3}.
\end{equation}

\end{cor}

The method of proof for Theorem \ref{mains} has its origin in \cite{AO}. Namely, the idea is
to first observe that since the functions ${}_{n+1}F_{n}(\lambda)$ are defined in terms of Jacobi sums, then one can express them as Gauss sums. One then applies the Gross-Koblitz formula \cite{gk} to express the Gauss sums in terms of $p$-adic Gamma functions. Using combinatorial properties of the $p$-adic Gamma function, Theorem \ref{mains} then follows. For an introduction to these methods, see \cite{cbms}. This general framework has been the basis for several recent results on supercongruences (see \cite{a}, \cite{kil}, \cite{lr}, \cite{mort}, \cite{mort1}, \cite{mort2}). Theorem \ref{mains} has recently been used to settle a conjecture of van Hamme (see \cite{mo}). Finally, the congruence in (\ref{leg}) appears to hold modulo $p^4$. This has been numerically confirmed for all primes less than $5000$.

The paper is organized as follows. In Section 2, we recall
properties of the $p$-adic Gamma function. In Section 3, we prove
Theorem \ref{mains}. In Section 4, we prove Corollary \ref{done}
using Theorem \ref{mains} and two non-trivial Harmonic sum
identities discovered using the computer summation
program~\SigmaP~\cite{Schneider:07a}. A description of the
non-trivial methods involved using the ~\SigmaP~ package is included
in Section 5. We should also mention that similar harmonic number
identities were discovered and proven in \cite{ps}. These types of
identities played an important role in the proof of Beukers'
supercongruence for Ap{\'e}ry numbers (see \cite{aeoz} or Theorem 7
in \cite{AO}).

\section{Preliminaries}

We first recall the definition of the $p$-adic Gamma function and list some of its main properties. For more details, see \cite{k}, \cite{murty}, or \cite{r}. Let $| \cdot |$ denote the $p$-adic absolute value on $\mathbb{Q}_{p}$. For $n \in \mathbb{N}$, we define

$$
\Gamma_{p}(n):= (-1)^n \prod_{\substack{j<n \\ (j,p)=1}} j.
$$

\noindent One can extend this function to all $x \in \mathbb{Z}_{p}$ upon setting

$$
\Gamma_{p}(x) := \lim_{n \to x} \Gamma_{p}(n).
$$

\noindent The following Proposition provides some of the main properties of $\Gamma_{p}$.

\begin{prop} \label{properties} Let $n \in \mathbb{N}$ and $x \in \mathbb{Z}_{p}$. Then
\begin{enumerate}

\item $\displaystyle \Gamma_{p}(0)=1$.

\item $\displaystyle \frac{\Gamma_{p}(x+1)}{\Gamma_{p}(x)} = \left \{ \begin{array}{l}
-x \quad \mbox{if $| x | =1$},\\
-1 \quad  \mbox{if $| x | < 1$.}
\end{array}
\right. \\$

\item If $0 \leq n \leq p-1$, then $n!=(-1)^{n+1} \Gamma_{p}(n+1)$.

\item $\displaystyle | \Gamma_{p}(x) | =1$.

\item Let $x_{0} \in [1, 2, \dotsc, p]$ be the constant term in the $p$-adic expansion of $x$. Then
$$\Gamma_{p}(x) \Gamma_{p}(1-x) = (-1)^{x_{0}}.$$

\item If $x \equiv y \pmod {p^n}$, then $\Gamma_{p}(x) \equiv \Gamma_{p}(y) \pmod {p^n}$.

\end{enumerate}

\end{prop}

For $x \in \mathbb{Z}_{p}$, we define

\begin{equation} \label{g1}
G_{1}(x):=\frac{\Gamma_{p}^{\prime} (x)}{\Gamma_{p}(x)}
\end{equation}

\noindent and

\begin{equation} \label{g2}
G_{2}(x):=\frac{\Gamma_{p}^{\prime \prime} (x)}{\Gamma_{p}(x)}.
\end{equation}

\noindent One can check that $G_{1}(x)$ and $G_{2}(x)$ are defined for all $x \in \mathbb{Z}_{p}$ using the fact that $\Gamma_{p}(x)$ is locally analytic and $| \Gamma_{p}(x) | =1$. We now mention some congruence properties of the $p$-adic Gamma function. For a proof of this result, see \cite{cde} or \cite{kil}.

\begin{prop} \label{properties1} Let $p \geq 7$ be prime, $x \in \mathbb{Z}_{p}$, and $z \in p\mathbb{Z}_{p}$. Then
\begin{enumerate}

\item $G_{1}(x)$, $G_{2}(x) \in \mathbb{Z}_{p}$.

\item We have

\begin{center}
$\displaystyle \Gamma_{p}(x + z) \equiv \Gamma_{p}(x) \Biggl ( 1 + zG_{1}(x) + \frac{z^2}{2} G_{2}(x) \Biggr ) \pmod {p^3}$.
\end{center}

\item $\Gamma_{p}^{\prime}(x+z) \equiv \Gamma_{p}^{\prime}(x) + z\Gamma_{p}^{\prime \prime}(x) \pmod  {p^2}$.

\end{enumerate}
\end{prop}

We also need the following combinatorial congruence which relates $\Gamma_{p}$ to certain binomial coefficients.

\begin{lem} \label{bc} If $p$ is an odd prime and $1 \leq j \leq \frac{p-1}{2}$, then
$$
-\phi_{p}(-1) (-1)^{j} \binom{\frac{p-1}{2} + j}{j}\binom{\frac{p-1}{2}}{j} \equiv \frac{\Gamma_{p}(\frac{1}{2} + j)^2}{\Gamma_{p}(1+j)^2} \pmod {p^2}.
$$
\end{lem}

\begin{proof}
By Proposition \ref{properties} (3) and (5), we have

$$
\begin{aligned}
-\phi_{p}(-1) (-1)^{j} \binom{\frac{p-1}{2} + j}{j}\binom{\frac{p-1}{2}}{j}  &=  -\phi_{p}(-1) (-1)^j \frac{(\frac{p-1}{2} + j)!}{j!^2 (\frac{p-1}{2} - j)!} \\
&=  \frac{\Gamma_{p}(\frac{1}{2} + j + \frac{p}{2}) \Gamma_{p}(\frac{1}{2} + j - \frac{p}{2})}{\Gamma_{p}(1+j)^2}. \\
 \end{aligned}
$$

\noindent Now, using Proposition \ref{properties1} (2), we have

$$
\begin{aligned}
& \displaystyle \Gamma_{p}\Bigl(\frac{1}{2} + j + \frac{p}{2}\Bigr) \Gamma_{p}\Bigl(\frac{1}{2} + j - \frac{p}{2}\Bigr) \\ & \\ & \equiv \Bigl\lbrace \Gamma_{p}\Bigl(\frac{1}{2} + j\Bigr) + \frac{p}{2}\Gamma_{p}^{\prime}\Bigl(\frac{1}{2} + j\Bigr)  \Bigr \rbrace  \Bigl\lbrace \Gamma_{p}\Bigl(\frac{1}{2} + j\Bigr) - \frac{p}{2} \Gamma_{p}^{\prime}\Bigl(\frac{1}{2} + j\Bigr)  \Bigr \rbrace \pmod {p^2} \\
& {} \\ & \equiv \Gamma_{p}\Bigl(\frac{1}{2} + j\Bigr)^2 \pmod{p^2} \\ \\
\end{aligned}
$$

\noindent and the result follows.

\end{proof}

\noindent Finally, we need to define

\begin{equation} \label{a}
A(j):= G_{1}(\tfrac{1}{2} + j) - G_{1}(1+j)
\end{equation}

\noindent and for a positive integer $n$

\begin{equation} \label{b}
\begin{aligned}
B(n, j):= & \tfrac{n+1}{2}\Big(G_{2}(\tfrac{1}{2} + j) - G_{2}(1 + j)\Big) + \tfrac{(n+1)n}{2}G_{1}(\tfrac{1}{2} + j)^2 \\
& + \tfrac{(n+1)(n+2)}{2}G_{1}(1+j)^2 - (n+1)^2 G_{1}(\tfrac{1}{2} + j)G_{1}(1+j).
\end{aligned}
\end{equation}

We require the following Lemma which relates $A(j)$ and $B(n, j)$ to generalized Harmonic sums. The proof is similar to Lemma 4.1 in \cite{kil} and thus is omitted.

\begin{lem} \label{har} Let $p$ be an odd prime and $0 \leq j \leq \frac{p-1}{2}$. Then

\begin{equation} \label{acon}
A(j) \equiv H_{\frac{p-1}{2} + j}^{(1)} - H_{j}^{(1)} + 2p \sum_{r=0}^{j-1} \frac{1}{(2r+1)^2} \pmod {p^2}
\end{equation}

\noindent and

\begin{equation} \label{bcon}
B(n, j) \equiv \tfrac{(n+1)^2}{2}\Bigl ( H_{\frac{p-1}{2} + j}^{(1)} - H_{j}^{(1)} \Bigr)^2 - (\tfrac{n+1}{2})\Bigl( H_{\frac{p-1}{2} + j}^{(2)} - H_{j}^{(2)} \Bigr) \pmod p.
\end{equation}

\end{lem}

\section{Proof of Theorem \ref{mains}}

We are now in a position to prove Theorem \ref{mains}.

\begin{proof}
Let $n \geq 3$ be odd. From  (\ref{def1}) and (\ref{def2}), we know that

$$
-p^{n} {}_{n+1}F_{n}(\lambda) =\frac{1}{1-p}\sum_{\chi} J(\phi, \chi)^{n+1} \bar{\chi}(\lambda)
$$

\noindent where $\bar{\chi}$ is the complex conjugate of $\chi$. After expressing the Jacobi sum $J(\phi, \chi)^{n+1}$ in terms of Gauss sums, we then apply the Gross-Koblitz formula \cite{gk} to get (see also \cite{AO} or \cite{lr})

\begin{equation} \label{nasty}
\begin{aligned}
-p^{n} {}_{n+1}F_{n}(\lambda) &= \frac{1}{1-p} \Bigg \{ \phi_{p}(\lambda) + (-\phi_{p}(-1))^{\tfrac{n+1}{2}} \Bigg( \sum_{j=0}^{\tfrac{p-3}{2}} \frac{\Gamma_{p}(\tfrac{j}{p-1})^{n+1}}{\Gamma_{p}(\tfrac{1}{2} + \tfrac{j}{p-1})^{n+1}} \omega^{j}(\lambda) \\
& + p^{n+1} \sum_{j=\tfrac{p+1}{2}}^{p-2}  \frac{\Gamma_{p}(\tfrac{j}{p-1})^{n+1}}{\Gamma_{p}(\tfrac{j}{p-1} -  \tfrac{1}{2})^{n+1}} \omega^{j}(\lambda) \Bigg) \Bigg \}.
\end{aligned}
\end{equation}

\noindent Here $\omega$ is the Teichm{\"u}ller character which satisfies
$\displaystyle \omega(\lambda) \equiv \lambda^{p^{s-1}} \pmod {p^s}$ and thus

$$\displaystyle \omega^{j}(\lambda) \equiv \lambda^{jp^{s-1}} \pmod {p^s}$$

\noindent  for $s \geq 1$. As $n \geq 3$ is odd, the second sum in (\ref{nasty}) vanishes modulo $p^3$.
As $\frac{1}{1-p} \equiv 1 + p + p^2 \pmod {p^3}$ and thus $\frac{j}{p-1} \equiv -j - jp - jp^2 \pmod {p^3}$, we apply parts (5) and (6) of Proposition \ref{properties} and reindex the summation to obtain

\begin{equation} \label{main}
\begin{aligned}
& -p^{n} {}_{n+1}F_{n}(\lambda)\\
& \equiv (1+ p + p^2) \Bigg \{ \phi_{p}(\lambda) + (-\phi_{p}(-1))^{\tfrac{n+1}{2}} \sum_{j=1}^{\tfrac{p-1}{2}} \frac{\Gamma_{p}(\tfrac{1}{2} + j + jp + jp^2)^{n+1}}{\Gamma_{p}(1 + j + jp + jp^2)^{n+1}} \omega^{\tfrac{p-1}{2} - j}(\lambda) \Bigg \}    \pmod {p^3} .
\end{aligned}
\end{equation}

\noindent By Proposition \ref{properties1} (2), we see that

$$
\begin{aligned}
\Gamma_{p}(x_{0} + j + jp + jp^2)^{n+1} & \equiv \Gamma_{p}(x_{0} + j)^{n+1} \Bigl[ 1+ (n+1)(jp + jp^2)G_{1}(x_{0} +j) \\
& + \tfrac{n+1}{2}(jp +jp^2)^2 \Bigl(G_{2}(x_{0} + j) + nG_{1}(x_{0} + j)^2\Bigr) \Bigr] \pmod {p^3}
\end{aligned}
$$

\noindent for $x_{0} \in \mathbb{Z}_{p}$. We expand the numerator and denominator of (\ref{main}) with $x_{0}=\frac{1}{2}$ and $x_{0}=1$ respectively. After multiplying the numerator and denominator by

$$
1 - (n+1)jpG_{1}(1+j) - \tfrac{n+1}{2} j^2p^2 \Bigl(G_{2}(1+j) - (n+2)G_{1}(1+j)^2 \Bigr) - (n+1) jp^2G_{1}(1+j),
$$

\noindent we get

\begin{equation} \label{main1}
\begin{aligned}
& -p^{n} {}_{n+1}F_{n}(\lambda) \\
& \equiv (1+ p + p^2) \Bigg \{ \phi_{p}(\lambda) +  (-\phi_{p}(-1))^{\tfrac{n+1}{2}} \sum_{j=1}^{\tfrac{p-1}{2}} \frac{\Gamma_{p}(\tfrac{1}{2} + j)^{n+1}}{\Gamma_{p}(1 + j)^{n+1}} \Bigl(1 + (n+1)jpA(j) \\
& + (n+1)jp^2 A(j) + j^2 p^2 B(n, j) \Bigr)  \omega^{\tfrac{p-1}{2} - j}(\lambda) \Bigg \}    \pmod {p^3}
\end{aligned}
\end{equation}

\noindent where $A(j)$ and $B(n, j)$ are defined by (\ref{a}) and (\ref{b}). We now need to consider the sums

\begin{equation} \label{coep2}
\phi_{p}(\lambda) + (-\phi_{p}(-1))^{\tfrac{n+1}{2}} \sum_{j=1}^{\tfrac{p-1}{2}} \frac{\Gamma_{p}(\tfrac{1}{2} + j)^{n+1}}{\Gamma_{p}(1 + j)^{n+1}} \Bigl(1 + 2(n+1)jA(j) + j^2 B(n, j) \Bigr)  \omega^{\tfrac{p-1}{2} - j}(\lambda),
\end{equation}

\begin{equation} \label{coep}
\phi_{p}(\lambda) + (-\phi_{p}(-1))^{\tfrac{n+1}{2}} \sum_{j=1}^{\tfrac{p-1}{2}} \frac{\Gamma_{p}(\tfrac{1}{2} + j)^{n+1}}{\Gamma_{p}(1 + j)^{n+1}} \Bigl(1 + (n+1)jA(j) \Bigr)  \omega^{\tfrac{p-1}{2} - j}(\lambda),
\end{equation}

\noindent and

\begin{equation} \label{coe1}
\phi_{p}(\lambda) + (-\phi_{p}(-1))^{\tfrac{n+1}{2}} \sum_{j=1}^{\tfrac{p-1}{2}} \frac{\Gamma_{p}(\tfrac{1}{2} + j)^{n+1}}{\Gamma_{p}(1 + j)^{n+1}}  \omega^{\tfrac{p-1}{2} - j}(\lambda)
\end{equation}

\noindent which are the coefficients of $p^2$, $p$, and $1$ respectively in (\ref{main1}).
Observe that as we want to determine ${}_{n+1}F_{n}(\lambda)$ mod $p^3$, it suffices to compute (\ref{coep2}) mod $p$, (\ref{coep}) mod $p^2$, and (\ref{coe1}) mod $p^3$. Also note that

\begin{equation} \label{tech}
\omega^{\tfrac{p-1}{2} - j}(\lambda)=\omega^{\tfrac{p-1}{2}}(\lambda)\omega^{-j}(\lambda)=\phi_{p}(\lambda) \omega^{-j}(\lambda)
\end{equation}

\noindent as $\omega$ is of order $p-1$. By Lemma 4.4 in \cite{kil}, we see that

\begin{equation} \label{yeah}
(\tfrac{n+1}{2}) j \Bigl(H_{\frac{p-1}{2} + j}^{(1)} - H_{\frac{p-1}{2} - j}^{(1)} \Bigr) \equiv -2(n+1)jp \sum_{r=0}^{j-1} \frac{1}{(2r+1)^2} \pmod {p^2}.
\end{equation}

\noindent By Lemma \ref{bc}, we have

\begin{equation} \label{mandy}
\Bigg [-\phi_{p}(-1) (-1)^{j} \binom{\frac{p-1}{2} + j}{j}\binom{\frac{p-1}{2}}{j} \Bigg ]^{\tfrac{n+1}{2}} \equiv \frac{\Gamma_{p}(\frac{1}{2} + j)^{n+1}}{\Gamma_{p}(1+j)^{n+1}} \pmod {p^2}
\end{equation}

\noindent and so after combining Lemma \ref{har}, (\ref{tech}), (\ref{yeah}), (\ref{mandy}) and accounting for $j=0$, then (\ref{coep2}) is congruent modulo $p$ to (\ref{x}) and (\ref{coep}) is
congruent modulo $p^2$ to (\ref{y}). Here we have used the fact that $\Gamma_{p}(1)^2 = 1$ and $\Gamma_{p}(\tfrac{1}{2})^2 = -\phi_{p}(-1)$
and thus for $n\geq 3$ odd

$$
\frac{\Gamma_{p}(\frac{1}{2})^{n+1}}{\Gamma_{p}(1)^{n+1}}=(-\phi_{p}(-1))^{\tfrac{n+1}{2}}.
$$

\noindent By Proposition 2.5 in \cite{mort}, we have

$$
\frac{\Gamma_{p}(\frac{1}{2} + j)^{2}}{\Gamma_{p}(1+j)^{2}}=-\phi_{p}(-1)\binom{2j}{j}^2 16^{-j}
$$

\noindent and thus using (\ref{tech}), we have that (\ref{coe1}) and (\ref{z}) are congruent modulo $p^3$, namely

\begin{equation} \label{equal}
\begin{aligned}
\phi_{p}(\lambda) + (-\phi_{p}(-1))^{\tfrac{n+1}{2}} & \sum_{j=1}^{\tfrac{p-1}{2}} \frac{\Gamma_{p}(\tfrac{1}{2} + j)^{n+1}}{\Gamma_{p}(1 + j)^{n+1}} \omega^{\tfrac{p-1}{2} - j}(\lambda) \\ & \equiv \phi_{p}(\lambda) \sum_{j=0}^{\tfrac{p-1}{2}} \binom{2j}{j}^{n+1} 16^{-j\big(\tfrac{n+1}{2}\big)} \lambda^{-jp^2} \pmod {p^3}.
\end{aligned}
\end{equation}

This proves the result for $n \geq 3$ odd. A similar argument applies to the case $n \geq 2$ is even. We now turn to the case $n=1$. By (\ref{nasty}), we need only consider the last sum

\begin{equation} \label{last}
-\phi_{p}(-1) p^2 \sum_{j=\tfrac{p+1}{2}}^{p-2}  \frac{\Gamma_{p}(\tfrac{j}{p-1})^2}{\Gamma_{p}(\tfrac{j}{p-1} - \tfrac{1}{2})^2} \omega^{j}(\lambda).
\end{equation}

\noindent By (5) and (6) of Proposition \ref{properties} and after reindexing the exponent of $\omega(\lambda)$, (\ref{last}) is equivalent modulo $p^3$ to

$$
 -\phi_{p}(-1)p^2 \sum_{j=\tfrac{p+1}{2}}^{p-2}  \frac{\Gamma_{p}(\tfrac{1}{2} + j)^2}{\Gamma_{p}(1+j)^2} (\tfrac{1}{2} + j)^2 \omega^{\tfrac{3p-3}{2} -j}(\lambda).
$$

\noindent By repeated use of Proposition \ref{properties}
(2) we have for $\frac{p+1}{2} \leq j \leq p-2$ that

$$
\begin{aligned}
& \Gamma_{p}\Bigl(\frac{1}{2} + \frac{p+1}{2}\Bigr)^2 \equiv 1 \pmod p \\
& \Gamma_{p}\Bigl(\frac{1}{2} + \frac{p+1}{2} + 1\Bigr)^2 \equiv \Bigl(\frac{p+2}{2}\Bigr)^2 \pmod p \\
& \Gamma_{p}\Bigl(\frac{1}{2} + \frac{p+1}{2} + 2\Bigr)^2 \equiv \Bigl(\frac{p+4}{2}\Bigr)^2 \Bigl(\frac{p+2}{2}\Bigr)^2 \pmod p \\
& \hspace{1.25in} \vdots \\
& \Gamma_{p}\Bigl(\frac{1}{2} + p-2\Bigr)^2 \equiv \Bigl(\frac{2p-5}{2}\Bigr)^2 \Bigl(\frac{2p-7}{2}\Bigr)^2 \cdots \Bigl(\frac{p+2}{2}\Bigr)^2 \pmod p. \\
\end{aligned}
$$

\noindent By Proposition \ref{properties} (3), $\Gamma_{p}(1+j)^2 = (j!)^2$. Also using the fact that $\lambda^{p-1} \equiv 1 \pmod p$, we have

$$
\begin{aligned}
\omega^{\tfrac{3p-3}{2}-j}(\lambda) & \equiv \lambda^{\tfrac{3p-3}{2} -j} \pmod p \\
& \equiv \lambda^{\tfrac{p-1}{2} + p-1-j} \pmod p \\
& \equiv \lambda^{\tfrac{p-1}{2} - j} \pmod p \\
\end{aligned}
$$

\noindent for $\frac{p+1}{2} \leq j \leq p-2$ and thus

$$
\begin{aligned}
& \sum_{j=\tfrac{p+1}{2}}^{p-2}  \frac{\Gamma_{p}(\tfrac{1}{2} + j)^2}{\Gamma_{p}(1+j)^2} (\tfrac{1}{2} + j)^2 \omega^{\tfrac{3p-3}{2} -j}(\lambda) \\
& \equiv \frac{1}{(\frac{p+1}{2})!^2} \Bigl(\frac{p}{2} + 1\Bigr)^2 \omega^{p-2}(\lambda) +  \dotso +  \frac{(\frac{p+2}{2})^2 \cdots (\frac{2p-5}{2})^2}{(p-2)!^2} \Bigl(\frac{1}{2} + p-2 \Bigr)^2 \omega^{\tfrac{p+1}{2}}(\lambda)   \pmod p\\
& \\
& \equiv \frac{1}{-\phi_{p}(-1)(\frac{1}{2})} (1)^2 \lambda^{-1} + \cdots + \frac{(1)^2 (2)^2 \cdots (\frac{p-5}{2})^2}{-\phi_{p}(-1) (p-2)^2 (p-3)^2 \cdots (\frac{p+1}{2})^2} \Bigl(\frac{p-3}{2}\Bigr)^2 \lambda^{-1 - \tfrac{p-5}{2}} \pmod p \\
& \\
& \equiv -\phi_{p}(-1) D(p, \lambda) \pmod p. \\
\end{aligned}
$$

\noindent This proves the result for $n=1$.
\end{proof}

\section{Proof of Corollary \ref{done}}

Theorem \ref{mains} can be used to obtain modulo $p^3$ supercongruences in various settings.
For example, Ap{\'e}ry numbers \cite{AO}, traces of Frobenius endomorphisms on elliptic curves \cite{koike}, \cite{ono}, and colored partition functions \cite{ono} all occur as special values of  ${}_{n+1}F_{n}(\lambda)$ for certain $n$ and $\lambda$. We do not mention these results here, choosing instead to illustrate with one example. We now prove Corollary \ref{done}.

\begin{proof} If $p$ is an odd prime, then by Section 3 in \cite{g1},

$$
p \cdot {}_{2}F_{1}(1)=-\Biggl(\frac{-1}{p}\Biggr)
$$

\noindent and so by (\ref{n1}),

$$
\Biggl(\frac{-1}{p}\Biggr) \equiv p^2 \Bigl[ X(p,1, 1) + D(p, 1) \Bigr] + pY(p,1, 1) + Z(p,1, 1) \pmod {p^3}.
$$

\noindent We now claim that

\begin{equation} \label{xd}
X(p, 1, 1) + D(p,1) + 1 \equiv 0 \pmod p.
\end{equation}

\noindent In order to verify (\ref{xd}), we first study $X(p,1,1)$. By (\ref{x}), we have

$$
\begin{aligned}
X(p, 1,1) & = \sum_{j=0}^{\tfrac{p-1}{2}} \binom{\frac{p-1}{2} +
j}{j} \binom{\frac{p-1}{2}}{j} (-1)^{j}
\Bigl(1 + 4j\Bigl(H_{\frac{p-1}{2} + j}^{(1)} - H_{j}^{(1)} \Bigr)  +\\
&\quad\quad\quad\quad j^2 \Bigl(2 \Bigl ( H_{\frac{p-1}{2} +
j}^{(1)} - H_{j}^{(1)} \Bigr)^2 - \Bigl( H_{\frac{p-1}{2} + j}^{(2)}
- H_{j}^{(2)} \Bigr) \Bigr)\Bigr).
\end{aligned}
$$

\noindent The identity

\begin{equation} \label{cool}
\sum_{k=0}^{n} (-1)^{k}\binom{n + k}{k} \binom{n}{k}\Bigl(1 +
2k\Bigl(H_{n+k} - H_{k} \Bigr) \Bigr)=(-1)^n(2n+1)
\end{equation}

\noindent was discovered using \SigmaP~ (see Lemma 2.2 in \cite{mort}). In particular, we find

\begin{equation*}
\sum_{j=0}^{\tfrac{p-1}{2}} \binom{\frac{p-1}{2} + j}{j}
\binom{\frac{p-1}{2}}{j} (-1)^{j} \Bigl(1 + 2j\Bigl(H_{\frac{p-1}{2}
+ j}^{(1)} - H_{j}^{(1)} \Bigr)  \Bigr) \equiv 0 \pmod p
\end{equation*}

\noindent and thus

\begin{equation} \label{yeah1}
\begin{aligned}
& X(p, 1,1) \\
& \equiv -1 - \sum_{j=1}^{\tfrac{p-1}{2}} \binom{\frac{p-1}{2} +
j}{j} \binom{\frac{p-1}{2}}{j} (-1)^{j} \\
& + \sum_{j=1}^{\tfrac{p-1}{2}} \binom{\frac{p-1}{2} +
j}{j} \binom{\frac{p-1}{2}}{j} (-1)^{j} \Bigl( j^2 \Bigl(2 \Bigl ( H_{\frac{p-1}{2} +
j}^{(1)} - H_{j}^{(1)} \Bigr)^2 - \Bigl( H_{\frac{p-1}{2} + j}^{(2)}
- H_{j}^{(2)} \Bigr) \Bigr)\Bigr) \pmod p. \\
\end{aligned}
\end{equation}

\noindent By Proposition \ref{properties1} and Lemma \ref{bc}, we also have

\begin{equation}
D(p,1) \equiv \sum_{j=1}^{\tfrac{p-3}{2}} \binom{\frac{p-1}{2} +
j}{j}^{-1} \binom{\frac{p-1}{2}}{j}^{-1} (-1)^{j} \pmod p.
\end{equation}

\noindent For positive integers $n$, the relation

\begin{equation} \label{new}
\begin{aligned}
& \sum_{j=1}^{n} \binom{n +
j}{j} \binom{n}{j} (-1)^{j} \Bigl( j^2 \Bigl(2 \Bigl ( H_{n+j}^{(1)} - H_{j}^{(1)} \Bigr)^2 - \Bigl( H_{n + j}^{(2)}
- H_{j}^{(2)} \Bigr) \Bigr)\Bigr) \\ & + \sum_{j=1}^{n} \binom{n +j}{j}^{-1} \binom{n}{j}^{-1} (-1)^{j} = n(-1 + 2n) (-1)^n
\end{aligned}
\end{equation}

\noindent was found using ~\SigmaP. Equation (\ref{xd}) now follows upon taking $n=\frac{p-1}{2}$ in
(\ref{new}), in the identity

\begin{equation} \label{old}
\sum_{j=1}^{n} \binom{n+j}{j}\binom{n}{j} (-1)^j = -1 + (-1)^n,
\end{equation}

\noindent and reducing modulo $p$. We now consider $Y(p, 1, 1)$. By (\ref{y}), we have

$$
\begin{aligned}
Y(p, 1, 1) = \sum_{j=0}^{\tfrac{p-1}{2}}  \binom{\frac{p-1}{2} + j}{j}\binom{\frac{p-1}{2}}{j} (-1)^{j} \Biggl(1 + j  &\Bigl(H_{\frac{p-1}{2} + j}^{(1)} - H_{j}^{(1)} \Bigr) \\
& +  j \Bigl(H_{\frac{p-1}{2} - j}^{(1)} - H_{j}^{(1)} \Bigr) \Biggr).
\end{aligned}
$$

\noindent For positive integers $n$, the relation

\begin{equation} \label{rel2}
\begin{aligned}
& \sum_{j=0}^{n} \binom{n + j}{j} \binom{n}{j} (-1)^{j} \Biggl(1+ j\Bigl (H_{n + j}^{(1)} + H_{n-j}^{(1)} -2H_j^{(1)}\Bigr )\Biggr) \\
& =(1+2n)\binom{2n}{n}(-1)^n - \frac{3}{2}n(1+n)(-1)^{n} \sum_{i=1}^{n} \frac{\binom{2i}{i}}{i} \\
 \end{aligned}
\end{equation}

\noindent was discovered using ~\SigmaP. Taking $n=\frac{p-1}{2}$ in (\ref{rel2}) and reducing mod $p^2$, we
have

\begin{equation} \label{yp}
Y(p, 1, 1) \equiv  p + \frac{3}{8} (-1)^{\tfrac{p-1}{2}} \sum_{i=1}^{\frac{p-1}{2}}  \binom{2i}{i} \frac{1}{i} \pmod {p^2}.
\end{equation}

\noindent Equation (\ref{leg}) then follows from  (\ref{z}), (\ref{xd}), and (\ref{yp}).
\end{proof}

\section{Finding and proving identities~\eqref{new}
and~\eqref{rel2} with ~\SigmaP}

An efficient algorithm to find and prove identities involving nested
definite and indefinite sum expressions, such as~\eqref{new}
and~\eqref{rel2}, has only recently been developed and implemented.
Inspired by hypergeometric summation~\cite{AequalB}, in particular
Zeilberger's creative telescoping method~\cite{Zeilberger:91}, and
Karr's indefinite summation algorithm~\cite{karr1,karr2} (extended
to definite summations), the second author developed and implemented
an algorithm using Mathematica to handle various summations. The
resulting package is called~\SigmaP. For a more detailed description
of the algorithms incorporated into~\SigmaP, please see
\cite{Schneider:07a}. Applications of this computer algebra package
include proving identities that arise in the enumeration of rhombus
tilings of a symmetric hexagon~\cite{fk,Schneider:04}, in the
verification of Stembridge's totally symmetric plane partitions
theorem~\cite{stem,aps}, and in certain Pad{\'e}
approximations~\cite{dpsw}. In this section, we illustrate how the
package can be used to discover and prove identities~\eqref{new}
and~\eqref{rel2}. For simplicity, we write $H_{k}$ for
$H_{k}^{(1)}$.

\subsection{Identity~\eqref{rel2}}

With \SigmaP\ we produce the following harmonic sum identities:

\begin{align}
\sum_{k=0}^n(-1)^k\binom{n+k}{k}\binom{n}{k}kH_{k}=&(-1)^nn(n+1)(2H_n-1),\label{Equ:Sumk}\\
\sum_{k=0}^n(-1)^k\binom{n+k}{k}\binom{n}{k}kH_{n+k}=&(-1)^n n(n+1)2H_n-(-1)^nn^2,\label{Equ:Sumn+k}\\
\intertext{and}
\sum_{k=0}^n(-1)^k\binom{n+k}{k}\binom{n}{k}kH_{n-k}=&-(-1)^n(n+1)^2+(-1)^n(2n+1)\binom{2n}{n}\label{Equ:Sumn-k}\\
&+2n(n+1)(-1)^nH_n-\frac{3}{2}n(n+1)(-1)^n\sum_{i=1}^n\frac{\binom{2i}{i}}{i}.\nonumber
\end{align}

\noindent Then, combining~\eqref{old}, \eqref{Equ:Sumn+k}
and~\eqref{Equ:Sumn-k} we arrive at identity~\eqref{rel2}.

\begin{rem}
Note that~\eqref{cool}, \eqref{Equ:Sumk} and \eqref{Equ:Sumn+k}
can be proved using hypergeometric techniques which appear
in \cite{Andrews74} and \cite{AU}. The key
observation 
is that differentiation of the rising factorial $(x)_k=x(x+1)\dots(x+k-1)$
(resp.\ $1/(x)_k$) in $x$ and afterwards substituting $x=1$
produces $(1)_kH_k$ (resp.\ $-H_k/(1)_k$).
With this fact, one can produce, e.g.,~\eqref{Equ:Sumk} by setting up
the identity
\begin{equation}\label{Equ:GaussApl}
\sum_{k=0}^n\frac{(-n)_k(n+1)_k}{k!(x)_k}k=
-\frac{n(n+1)}{x}{}_{2}F_{1}(1-n,n+2;x+1;1)
=-\frac{n(n+1)}{x}\frac{(x-n-1)_{n-1}}{(x+1)_{n-1}}
\end{equation}
with Gauss' theorem, differentiating~\eqref{Equ:GaussApl} in $x$, and setting $x=1$.
Obviously, the successful  application of this technique relies
on the fact that one knows the underlying hypergeometric identity
such as~\eqref{Equ:GaussApl} for the particular case~\eqref{Equ:Sumk}. It would be interesting to see
proofs of identity~\eqref{Equ:Sumn-k}, in particular
of identity~\eqref{new}, along the lines sketched above.
Recently, a skillful application of partial fraction decomposition has been
used in~\cite{Prod:07} to derive identity~\eqref{Equ:Sumn-k}, but so far
no proof of identity~\eqref{new} has been found.
\end{rem}

Subsequently, we illustrate the computation steps for identity~\eqref{Equ:Sumn-k}
which can be executed in a straightforward manner.
After loading the package

\begin{mma}
\In << |Sigma.m| \\
\Print Sigma - A summation package by Carsten Schneider
\copyright\ RISC-Linz\\
\end{mma}

\noindent into the computer algebra system Mathematica, we insert
the sum in question:

\begin{mma}
\In S=SigmaSum[SigmaPower[-1, k]k SigmaBinomial[n+k,k]\newline
\hspace*{2cm}SigmaBinomial[n,k]SigmaHNumber[n - k], \{k, 0, n\}]\\
\Out \sum_{k=0}^n(-1)^kk\binom{n+k}{k}\binom{n}{k}H_{n-k}\\
\end{mma}

\begin{rem} Various functions support the user, like {\tt
SigmaSum} for sums, {\tt SigmaPower} for powers, {\tt SigmaBinomial}
for binomials, or {\tt SigmaHNumber} for harmonic numbers.
\end{rem}

\smallskip

Next, we compute a recurrence relation for the given sum ${\tt S}$
by inputting:

\begin{mma}\MLabel{Out:RecS1}
\In rec=GenerateRecurrence[{\tt S}]\\
\Out (n+2)(2n+1)(n+1)^2
SUM[n+2]+2(n+3)(2n^2+4n+1)(n+1)SUM[n+1]+\newline
\hspace*{0.41cm}(n+1)(n+2)(n+3)(2n+3)SUM[n]==-(2n+1)(2n+3)(3n+1)(3n+4)(-1)^n\binom{2n}{n}\\
\end{mma}

\noindent This means that $\text{SUM[n]}(={\tt S}= \displaystyle
\sum_{k=0}^n(-1)^kk\tbinom{n+k}{k}\tbinom{n}{k}H_{n-k})$
satisfies~\myOut{\ref{Out:RecS1}}.

\medskip

\noindent{\it Proof of~\myOut{\ref{Out:RecS1}}:} Define
$f(n,k):=(-1)^kk\binom{n+k}{k}\binom{n}{k}H_{n-k}$. The correctness
follows by the creative telescoping equation
\begin{equation}\label{Equ:CreaTeleEqu}
g(n,k+1)-g(n,k)=c_0(n)f(n,k)+c_1(n)f(n+1,k)+c_2(n)f(n+2,k)
\end{equation}
and the proof certificate $c_0(n)=(n+2)(n+3)(2n+3)$,
$c_1(n)=2(n+3)\left(2n^2+4n+1\right)$, $c_2(n)=(n+1)(n+2)(2n+1)$
and\small
\begin{multline*}
g(n,k)=(k-1)k^2\Big(2H_{n-k}(k-n-2)(k-n-1)(n+1)(k(4n+7)-2(2
n^3+10n^2+17n+10))+\\
(-k-n-1)(16n^4+88n^3+179n^2+163n+2k^2(4n^2+11n+7)-k(24n^3+\\
\hspace*{2cm}98n^2+131n+59)+58)\Big)
(-1)^k\binom{n+k}{k}\binom{n}{k}\Big/\big((n+1)(-k+n+1)^2(-k+n+2)^2\big)
\end{multline*}
\normalsize delivered by~\SigmaP. We verify~\eqref{Equ:CreaTeleEqu} as follows. Express
$g(n,k+1)$ in terms of $h(n,k)=(-1)^k\binom{n+k}{k}\binom{n}{k}$ and
$H_{n-k}$ by using the relations
$$h(n,k+1)=-\frac{(n-k)(n+k+1)}{(k+1)^2}h(n,k)$$
\noindent and
$$H_{n-k-1}=H_{n-k}-\frac{1}{n-k}.$$
\noindent Similarly, express $f(n+i,k)$ in
terms of $h(n,k)$ and $H_{n-k}$ by using the relations
$$h(n+1,k)=\frac{n+k+1}{n-k+1}h(n,k)$$
\noindent and
$$H_{n-k+1}=H_{n-k}+\frac{1}{n-k+1}.$$
\noindent Then~\eqref{Equ:CreaTeleEqu}
can be checked directly.
Summing~\eqref{Equ:CreaTeleEqu} over $k$ from $0$ to $n$
produces~\myOut{\ref{Out:RecS1}}.\qed

\medskip

Next, we solve the recurrence relation~\myOut{\ref{Out:RecS1}} by
typing in:

\begin{mma}
\In recSol=SolveRecurrence[rec[[1]], SUM[n]]\\
\Out\{\{0,n(1+n)(-1)^n\},\{0,(1+n)(-1)^n\big(-1+2n\sum_{i=1}^n\frac{1}{i}\big)\},\newline
\hspace*{1cm}\{1,-\frac{1}{2}(-1)^n\big(-2(1+2n)\binom{2n}{n}+3n(1+n)\sum_{i=1}^n\frac{\binom{2i}{i}}{i}\big)\}\}\\
\end{mma}

\noindent The result has to be interpreted as follows. \SigmaP\
finds two linearly independent solutions $h_1(n)=n(1+n)(-1)^n$ and
$h_2(n)=\displaystyle (1+n)(-1)^n\big(-1+2n\sum_{i=1}^n\frac{1}{i}\big)$ of the
the homogeneous version of~\myOut{\ref{Out:RecS1}} (indicated by the
$0$ in front) plus one particular solution
$$p(n)=-\frac{1}{2}(-1)^n\big(-2(1+2n)\binom{2n}{n}+3n(1+n)\sum_{i=1}^n\frac{\binom{2i}{i}}{i}\big)$$
of the input recurrence itself (indicated by the $1$ in front). The
correctness of the result can be easily verified by using, e.g., the
relation

$$\sum_{i=1}^{n+1}\binom{2i}{i} \frac{1}{i} =\sum_{i=1}^n\binom{2i}{i} \frac{1}{i} +
\frac{2(2n+1)}{(n+1)^2}\binom{2n}{n}.$$

\smallskip

Finally, by taking all linear combinations
$c_1h_1(n)+c_2h_2(n)+p(n)$ for constants $c_1$ and $c_2$, free of
$n$, we obtain all solutions of~\myOut{\ref{Out:RecS1}}. Hence, by
considering the first two initial values of~{\tt S} we can discover
and prove~\eqref{Equ:Sumn-k}:

\begin{mma}
\In FindLinearCombination[recSol, {\tt S}, 2]\\
\Out
-(-1)^n(1+n)^2+(-1)^n(1+2n)\binom{2n}{n}+2n(1+n)(-1)^n\sum_{i=1}^n\frac{1}{i}
-\frac{3}{2}n(1+n)(-1)^n\sum_{i=1}^n\frac{\binom{2i}{i}}{i}\\
\end{mma}

\begin{rem} Looking at the
identities~\eqref{old},\eqref{Equ:Sumn+k} and~\eqref{Equ:Sumk} one
immediately sees that the
combination~$\eqref{old}+2\cdot\eqref{Equ:Sumn+k}-2\cdot\eqref{Equ:Sumk}$
produces ~\eqref{cool}. Since the sums can be combined so nicely, we
had also the hope to find a solution for the sum
$$S_{\lambda}(n):=\sum_{k=0}^{n} (-\lambda)^{k}\binom{n + k}{k} \binom{n}{k}\Bigl(1 +
2k\Bigl(H_{n+k} - H_{k} \Bigr) \Bigr).$$ \SigmaP\ was able to compute the
recurrence
\begin{multline*}
(n+2)^2S_{\lambda}(n)+(2\lambda-1)\left(4n^2+18n+21\right)S_{\lambda}(n+1)\\
+\left(16n^2\lambda^2+80n\lambda^2+100\lambda^2-16n^2\lambda-80n\lambda-100\lambda+6
n^2+30n+39\right)S_{\lambda}(n+2)\\
+(2\lambda-1)\left(4n^2+22n+31\right)S_{\lambda}(n+3)+(n+3)^2S_{\lambda}(n+4)=0,
\end{multline*}
but failed to find any solution for a generic value $\lambda$.
Interesting enough, choosing $\lambda=\frac{1}{2}$ the recurrences
gets much simpler. In particular, this indicates that considering
the sums $S_{\frac{1}{2}}(2n)$ and $S_{\frac{1}{2}}(2n+1)$
separately, one can compute recurrences of order $2$ for each of
them. Indeed, applying the mechanism from above for each of the sums
gives (two different) recurrences of order two. Luckily, we can even
solve the recurrences which yields
\begin{align*}
S_{\frac{1}{2}}(2n)&=\frac{(-1)^n 2^{2 n} (n!)^2}{(2 n)!}\\
\intertext{and} S_{\frac{1}{2}}(2n+1)&=\frac{(-1)^n(2 n)!}{2^{2
n}(n!)^2}\Big((2 n+1)\big(H_n-H_{2n}\big)-1\Big).
\end{align*}

\end{rem}

\subsection{Identity~\eqref{new}}

Finally, we derive the two identities
\begin{equation} \label{Equ:AlgSum1}
\begin{aligned}
\sum_{k=1}^{n}& \binom{n + k}{k}\binom{n}{k} (-1)^{k}\Biggl(2k^2 \Bigl ( H_{n + k} - H_{k} \Bigr)^2 - k^2 \Bigl( H_{n + k}^{(2)} - H_{k}^{(2)} \Bigr) \Biggr) \\
=&(1+n)^2 (-2-2n+n^2)\frac{(n!)^2(-1)^n}{2(2+2n)!}
+ \frac{1}{4}n(-4 + 11n + 6n^2 + 3n^3)(-1)^n\\
& -\frac{1}{2}(-1+n +n^2) + \frac{3}{2}n^2 (1+n)^2(-1)^n
\sum_{i=1}^{n} \frac{i!^2}{(2+2i)!}
+ n^2 (1+n)^2(-1)^n \sum_{i=1}^{n} \frac{(-1)^i}{i^2} \\
\end{aligned}
\end{equation}
and
\begin{equation} \label{Equ:AlgSum2}
\begin{aligned}
\sum_{k=1}^{n}& \binom{n + k}{k}^{-1}\binom{n}{k}^{-1} (-1)^{k}=n(2n-1)(-1)^n\\
&-\Big((1+n)^2 (-2-2n+n^2)\frac{(n!)^2(-1)^n}{2(2+2n)!}
+ \frac{1}{4}n(-4 + 11n + 6n^2 + 3n^3)(-1)^n\\
& -\frac{1}{2}(-1+n +n^2) + \frac{3}{2}n^2 (1+n)^2(-1)^n
\sum_{i=1}^{n} \frac{i!^2}{(2+2i)!}
+ n^2 (1+n)^2(-1)^n \sum_{i=1}^{n} \frac{(-1)^i}{i^2}\Big)\\
\end{aligned}
\end{equation}
which immediately gives identity~\eqref{new}.

One option is to follow the same strategy as above: We can compute a
recurrence for
\begin{mma}
\In mySum=\sum_{k=0}^{n}(-1)^{k}\binom{n+k}{k}\binom{n}{k}
\Biggl( 2k^2\Bigl(H_{n+k}-H_{k}\Bigr)^2-k^2\Bigl(H_{n+k}^{(2)}-H_{k}^{(2)}\Bigr)\Biggr)\\
\end{mma}
\noindent and can solve the derived recurrence to find the right
hand side of~\eqref{Equ:AlgSum1}. But, since the found recurrence
relation is rather big (it has order 4), and the proof certificate
is even bigger (it fills about one page), we follow a refined
strategy presented in~\cite{ps} and~\cite{Schneider:07b}. Namely,
by running our creative telescoping algorithm with the additional
option {$\tt SimplifyByExt\rightarrow DepthNumber$} we can find a
recurrence of smaller order (order one!):

\begin{mma}\MLabel{MMA:AlgebraicSumRec}
\In rec = GenerateRecurrence[mySum, SimplifyByExt\rightarrow DepthNumber]\\
\Out 2(2n+1)(n+2)^2 SUM[n]+2(2n+1)n^2SUM[n+1]==\newline
\hspace*{-0.5cm}4(1+2n)+n^2(n+1)(n+2)(3n+2)
\sum_{i=0}^n\frac{(-1)^i\tbinom{n+i}{i}\tbinom{n}{i}}{(n+i)^2}+
2n\left(4n^2+3n-4\right)(2n+1)
\sum_{i=0}^n(-1)^i\tbinom{n+i}{i}\tbinom{n}{i}+\newline
\hspace*{1cm}8(n-1)n(n+1)(n+2)(2n+1)
\Big(\sum_{i=0}^n(-1)^i\tbinom{n+i}{i}\tbinom{n}{i}H_{n+i}-\sum_{i=0}^n(-1)^i\tbinom{n+i}{i}\tbinom{n}{i}H_{i}\Big)\\
\end{mma}

\medskip

\noindent{\it Proof of~\myOut{\ref{MMA:AlgebraicSumRec}}:} Define
$f(n,k)=(-1)^{k}\binom{n+k}{k}\binom{n}{k}\big(
2k^2\big(H_{n+k}-H_{k}\big)^2-k^2\big(H_{n+k}^{(2)}-H_{k}^{(2)}\big)\big)$.
Then the correctness of~\myOut{\ref{MMA:AlgebraicSumRec}} follows by
the creative telescoping equation
\begin{equation}\label{Equ:CreaTeleOrder1}
g(n,k+1)-g(n,k)=c_0(n)f(n,k)+c_1(n)f(n+1,k)
\end{equation}
with the proof certificate $c_0(n)=2(n+2)^2(2n+1)$,
$c_1(n)=2n^2(2n+1)$ and\footnotesize
\begin{multline*}
g(n,k)=\Bigg(4(k-1)^2n(n+1)^2(2 n+1)k^2\bigg(2H_k^2-4H_{n+k}H_k+2
H_{n+k}^2+H_k^{(2)}-H^{(2)}_{n+k}\bigg)-\\
\Big(n(n+2)k^3-(n^3+2n^2+2n+2)k^2-
(n+1)^2(n^2-2)k+n(n+1)^2(n^2+n-2)\Big)8n(n+1)\\(2n+1)
\bigg(H_k-H_{n+k}\bigg)+(16n^5+48n^4+29n^3+ 14n^2+20n+8)k^2+n(n+1)^2(16n^4+23n^3+\\
n^2+12n+8)-(32n^6+101n^5+98n^4+55n^3+
54n^2+36n+8)k\Bigg)\frac{(-1)^k\binom{n+k}{k}\binom{n}{k}}{(1-k+n)n(1+n)}+\\
n^2(n+1)(n+2)(3n+2)
\sum_{i=0}^k\tfrac{(-1)^i\tbinom{n+i}{i}\tbinom{n}{i}}{(n+i)^2}+
2n\left(4n^2+3n-4\right)(2n+1)
\sum_{i=0}^k(-1)^i\tbinom{n+i}{i}\tbinom{n}{i}\\
+8(n-1)n(n+1)(n+2)(2n+1)
\Big(\sum_{i=0}^k(-1)^i\tbinom{n+i}{i}\tbinom{n}{i}H_{n+i}-\sum_{i=0}^k(-1)^i\tbinom{n+i}{i}\tbinom{n}{i}H_{i}\Big).
\end{multline*}
\normalsize Since the sums and products inside of $g(n,k)$ are all
indefinite, e.g., we can apply the relation

$$\sum_{i=0}^{k+1}(-1)^i\tbinom{n+i}{i}\tbinom{n}{i}=\sum_{i=0}^{k}(-1)^i\tbinom{n+i}{i}\tbinom{n}{i}-
(-1)^k\tfrac{(k-n) (k+n+1)}{(k+1)^2}\tbinom{n+k}{k}\tbinom{n}{k},$$

\noindent the verification of~\eqref{Equ:CreaTeleOrder1} is
immediate. Summing~\eqref{Equ:CreaTeleOrder1} over $k$ from $0$ to
$n$ produces~\myOut{\ref{MMA:AlgebraicSumRec}}.\qed

\medskip

At first glance the recurrence~\myOut{\ref{MMA:AlgebraicSumRec}}
seems to be disappointing: we start with the definite sum {\tt
mySum}, and end up with a recurrence again involving definite sums.
But, these sums are much simpler than the input sum. In particular,
facilitating again \SigmaP, we can produce mechanically the
identities
$$\sum_{i=0}^n(-1)^i\frac{\tbinom{n+i}{i}\tbinom{n}{i}}{(n+i)^2}=-(-1)^n\frac{n!^2}{n^2(2n)!}$$
and
$$\sum_{k=0}^n(-1)^k\binom{n+k}{k}\binom{n}{k}H_{k}=\sum_{k=0}^n(-1)^k\binom{n+k}{k}\binom{n}{k}H_{n+k}=(-1)^n2H_n.$$
Using in addition~\eqref{old}, we can simplify the recurrence to
\begin{mma}
\In
rec=rec/.\{\sum_{i=0}^n(-1)^i\tfrac{\tbinom{n+i}{i}\tbinom{n}{i}}{(n+i)^2}\to-(-1)^n\tfrac{n!^2}{n^2(2n)!},
\sum_{i=0}^n(-1)^i\tbinom{n+i}{i}\tbinom{n}{i}\to(-1)^n,\newline
\hspace*{2cm}\sum_{i=0}^n(-1)^i\tbinom{n+i}{i}\tbinom{n}{i}H_{i}\to\sum_{i=0}^n(-1)^i\tbinom{n+i}{i}\tbinom{n}{i}H_{n+i}\}\\
\Out 2(2n+1)n^2\text{SUM}[n+1]+2(n+2)^2(2n+1)\text{SUM}[n]==\newline
\hspace*{1cm}4(2n+1)-\frac{(-1)^n(n+1)(n+2)(3n+2)(n!)^2}{(2n)!}+2(-1)^nn(2n+1)\left(4n^2+3n-4\right).\\
\end{mma}

\noindent Given this recurrence, one can directly read off its
solution. With some simplifications \SigmaP\ yields:

\begin{mma}
\In
recSol=SolveRecurrence[rec[[1]],SUM[n],SimpleSumRepresentation\to
True]\\
\Out \{\{0,(-1)^n n^2 (n+1)^2, \{1,(1+n)^2
(-2-2n+n^2)\frac{(n!)^2(-1)^n}{2(2+2n)!}-\frac{1}{2}(-1+n
+n^2)+\newline \frac{1}{4}n(-4 + 11n + 6n^2 +
3n^3)(-1)^n+\frac{3}{2}n^2 (1+n)^2(-1)^n \sum_{i=1}^{n}
\frac{i!^2}{(2+2i)!}
+ n^2 (1+n)^2(-1)^n \sum_{i=1}^{n} \frac{(-1)^i}{i^2}\}\}\\
\end{mma}

\noindent Looking at the first initial values we end up at the
identity~\eqref{Equ:AlgSum1}.

Finally, we attack the sum $S(n)$ on the left hand side
of~\eqref{Equ:AlgSum2}. Namely, we derive the recurrence
$$(n+2)S(n)+n^2S(n)=\frac{(-1)^n
(n+1)^2 (n+2) (3 n+2)
   (n!)^2}{(2 n+2)!}-2;$$
we remark that for this hypergeometric sum any implementation of
Zeilberger's algorithm~\cite{Zeilberger:91} could do the job. Using
this information, we find as above, the right hand side
of~\eqref{Equ:AlgSum2}.

\section*{Acknowledgments}
The first author would like to thank the Institut des Hautes {\'E}tudes
Scientifiques for their hospitality and support during the preparation
of this paper. The authors also thank Ken Ono for his comments on a preliminary version of the paper, Dermot McCarthy for his careful reading, and the referee for helpful suggestions.

\end{document}